\documentclass[a4paper,10pt]{article}
\usepackage{amsmath,amsthm,amssymb}
\usepackage{graphicx,subfigure}
\usepackage{url}

\usepackage[colorlinks,bookmarksopen,bookmarksnumbered,citecolor=red,urlcolor=red]{hyperref}
\hypersetup{pdftitle={Parallelization of a relaxation scheme modelling the bedload transport of sediments in shallow water flow},
bookmarks=true,
pdftoolbar=true,
pdfmenubar=true,
pdfauthor={E. Audusse, O.Delestre, M.-H. Le, M. Masson-Fauchier, P. Navaro, R. Serra},
pdfsubject={applied mathematics, hydraulics, erosion},
pdfcreator={Delestre},
pdfproducer={Delestre},
pdfkeywords={Shallow-Water equation} {Saint-Venant} {Saint-Venant-Exner} {relaxation method} {dune}
 {anti-dune} {antidune} {MPI} {parallelization} {Grass} {Meyer-Peter & M\"uller} {Strickler} {Darcy-Weisbach}
 {Validation of numerical method}}

\title{Parallelization of a relaxation scheme modelling the bedload transport of sediments in shallow water flow}

\author{{E. Audusse}\footnote{BANG Project, INRIA-Paris-Rocquencourt \& LAGA, Universiy Paris 13 Nord, France, e-mail :
 audusse@math.univ-paris13.fr}, O. Delestre\footnote{Lab. J.A. Dieudonn\'e \& EPU Nice Sophia, University of
 Nice, France, e-mail : delestre@math.unice.fr}\;,
 {M.-H. Le}\footnote{MAPMO, University of Orl\'eans, France, e-mail : Minh.Hoang.Le@math.cnrs.fr}\,,
 M. Masson-Fauchier\footnote{EPU Nice Sophia, University of Nice, France},\\
 P. Navaro\footnote{IRMA, UMR 7501 CNRS, Unistra, France, e-mail : navaro@math.unistra.fr}
 and R. Serra\footnote{EPU Nice Sophia, University of Nice, France}}

\begin{document}
\maketitle

\begin{abstract} 
In this work we are interested in numerical simulations for bedload erosion processes. We present a relaxation solver that we apply to moving
 dunes test cases in one and two dimensions. In particular we retrieve the so-called anti-dune process that is well described in the experiments.
 In order to be able to run 2D test cases with reasonable CPU time, we also describe and apply a parallelization procedure by using domain
 decomposition based on the classical MPI library.

Nous nous int\'eressons dans ce travail \`a la simulation num\'erique des processus d'\'erosion par charriage en rivi\`ere. Nous pr\'esentons
 un sch\'ema num\'erique bas\'e sur un mod\`ele de relaxation et nous l'appliquons \`a des probl\`emes de d\'eplacement de dunes en une et deux
 dimensions d'espace. Nous \'etudions en particulier les ph\'enom\`enes d'antidune pr\'esent\'es dans la litt\'erature hydros\'edimentaire.
 Nous pr\'esentons enfin une proc\'edure de parall\'elisation bas\'ee sur l'utilisation de la librairie MPI et qui nous permet de simuler
 des processus bidimensionnels avec des temps CPU raisonnables.
\end{abstract}

\section{Introduction}
Soil erosion is a complex phenomenon affected by many factors such as climate, topography, soil characteristics, vegetation and anthropogenic
 activities such as cultivation practices. Erosion process can be described in three stages: detachment, transport and deposition. The detachment
 occurs when the flow shear stress or the kinetic energy of raindrop exceeds the cohesive strength of the soil particles. Once detached,
 the sediments can be transported downstream as non-cohesive sediment before its deposition.

The movement of sediments occurs in two main modes called bedload and suspended load. The bedload particles are located in a few grain
 diameters thick layer situated on the soil. The velocities of these particles are less than the flow velocity. At the opposite, the suspended
 particles are transported in the flow without contact with the bed. Sediments finer than $0.2~{\rm mm}$ which are transported in suspension
 are rarely included/considered in bedload. The distinction between these two modes of sediment transport is blurred because they occur together. 

As soil erosion by water continues to be a serious problem throughout the world, the development of improved soil erosion prediction technology
 is required. With the increase of computing powers in the last years, there has been a rapid increase in the erosion and sediment transport
 simulations through the use of computer models. The models describing erosion process are available at different level of complexity. In
 general, erosion process is described by the equations of evolution based on the principle of conservation. These equations are derived at
 small scale under physical assumptions.

\section{Bedload modelling}

In this paper, we focus on the modelling of the morphodynamic process where the solid transport is only characterized by bedload and thus the
 suspended load has been ignored. The governing equations are often given by coupling the shallow water equations, describing the flow
 routing \cite{Gerbeau01}, with the Exner equation \cite{Exner1925} expressing the mass conservation of sediment layer. The one-dimensional
 system may be written in the form
\begin{equation}\label{bedload system}
 \left\{
\begin{array}{l}
\partial_t h +\partial_x (hu) = 0,\\
\partial_t (hu) + \partial_x(hu^2 + gh^2/2) = -gh\partial_x z_b,\\
\partial_t z_b + \partial_x q_b = 0, 
\end{array}
\right.
\end{equation}
where $h$ is the water depth, $u$ the flow velocity, $z_b$ the thickness of sediment layer, $q_b$ the volumetric bedload sediment transport rate
 per unit time and width and $g$ the acceleration due to gravity.

The expression of bedload flux $q_b$ is necessary to close the model~\eqref{bedload system}. Many researches have developed different empirical
 formul\ae{} to predict and to estimate $q_b$. One of the simplest expression was proposed by Grass \cite{grass81sediment} where $q_b$ is a
 function of the flow velocity and a dimensional $A_g$ constant, called interaction constant, that encompasses the effects of grain size and
 kinematic viscosity and is usually determined from experimental data:
\begin{equation}\label{bedload Grass}
q_b = A_gu|u|^{m_g-1},\quad 1\leq m_g\leq 4.
\end{equation}
The usual value of the exponent $m_g$ is set to $m_g = 3$. If $A_g = 0$ then we have a solid bed (no sediment transport) and we recover the
 standard shallow water equations. When $A_g$ is near zero, there is a small interaction between the fluid and the bed, while if $A_g$ is near
 one the interaction is larger. Note that one of the main characteristics of this model is that threshold value to initiate the motion of
 sediment is set to zero, so the sediment transport begins at the same instant that beginning fluid motion.

In practice, $q_b$ is usually represented under the non-dimensional form $q^*_b$ as a function of the dimensionless shear stress $\tau^*_b$ and a
 threshold value $\tau^*_{cr}$, {\it i.e.}
\begin{equation*}\label{chap2 non dimensional parameter}
q_b = q^*_b\sqrt{(s-1)gd_s^3},\quad q^*_b = q^*_b(\tau^*_b,\tau^*_{cr}),\quad \tau^*_b =\frac{\tau_b}{(\rho_s - \rho) gd_s},
\end{equation*}
where $\tau_b=\rho ghS_f$ is the bottom shear stress, {\it i.e.} the force of water acting on the bed during its routing, $s = \rho_s/\rho$ the
 density relative of sediment in water and $d_s$ the average value of sediment diameters. The friction term $S_f$ can be quantified by different
 empirical laws such as the Darcy-Weisbach or Manning formul\ae, {\it i.e.}
\begin{itemize}
 \item[$\bullet$] Darcy-Weisbach: $S_f =\frac{fu|u|}{8gh}$,
 \item[$\bullet$] Manning: $S_f = \frac{n^2u|u|}{h^{4/3}}$,
\end{itemize}
where $f$ and $n$ are the Darcy-Weisbach and the Manning coefficients respectively. The threshold value $\tau^*_{cr}$ depends on the physical
 properties of sediment and is usually computed experimentally. One of the first works on this topic was done by Shields \cite{Shields1936}.
  
\begin{figure}[ht!]
\centering
\includegraphics[width=0.85\linewidth]{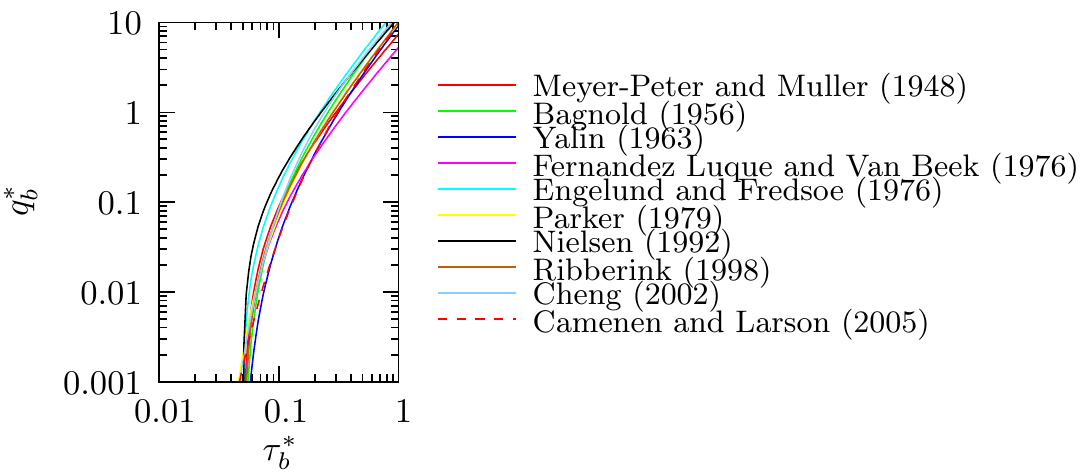}
\caption{Plot of several bedload functions found in the literature \cite{Garcia2008}, $\tau^*_{cr} = 0.05$.}
\label{bedload formula}
\end{figure}

The followings expressions, illustrated by Fig.~\ref{bedload formula}, have been often applied
 \cite{meyer-peter48formulas,Luque1976,nielsen92coastal,Camenen2005}:
\begin{align}
&\text{Meyer-Peter \& M\"uller (1948): } q^*_b =8(\tau^*_b-\tau^*_{cr})_+^{3/2} \label{MPM}\\
&\text{Fern\'andez Luque \& Van Beek (1976): } q^*_b =5.7(\tau^*_b-\tau^*_{cr})_+^{3/2} \label{FLVB}\\
&\text{Nielsen (1992): } q^*_b =12\sqrt{\tau^*_b}(\tau^*_b-\tau^*_{cr})_+\label{Nielsen}\\
&\text{Ribberink (1998): } q^*_b =11(\tau^*_b-\tau^*_{cr})_+^{1.65}\label{Ribberink}\\
&\text{Camenen and Larson (2005): } q^*_b=12(\tau^*_b)^{1.5}\exp \left(-4.5\tau^*_{cr}/\tau^*_b\right)\label{Camenen}
\end{align}

Upon initiation of flow routing, the topography $z_b$ becomes unstable resulting the fluid-sediment interaction. Various bedforms can be
 occurred, in particular the formation of {\it ripples, dunes} and {\it anti-dunes} (see e.g.~\cite{Kennedy1963,Kennedy1969}). Ripples are
 distinguished from dunes by their much smaller scale. Indeed, ripples and dunes often co-exist, with ripples forming on the larger dunes.
 At very low flow rates, ripples form on the bed, and as the flow rate increases, these are replaced by the longer wavelength and larger
 amplitude dunes which migrate slowly downstream. Ripples and dunes are observed in fluvial regime (the Froude number $Fr < 1$). When
 the flow becomes supercritical ($Fr>1$), dunes migrate upstream and called anti-dunes.

Let us show quickly how dunes and anti-dunes can be reproduced by the non-linear coupled system~\eqref{bedload system} (see \cite{Cordier2011}
 for more details). Rewriting the system in quasi-linear form as
\begin{equation*}
\partial_t W+A(W)\partial_x W= 0
\end{equation*}
where $W = (h, hu,z_b)$ and $A(W)$ the matrix of transport coefficients being
\begin{equation*}
A(W) =
\begin{pmatrix}
0 & 1 & 0\\
gh-u^2& 2u & gh\\
\alpha & \beta & 0
\end{pmatrix}
\end{equation*}
with $\alpha=\dfrac{\partial q_b}{\partial h}$ and $\beta=\dfrac{\partial q_b}{\partial hu}$. In~\cite{Cordier2011}, we can find two type of
 relations between $\alpha$ and $\beta$ for different formulas of $q_b$ (equations \eqref{bedload Grass}-\eqref{Camenen}) which depend on
 the friction laws:
\begin{itemize}
 \item with Darcy-Weisbach's law: $\alpha=-u\beta$,
 \item with Manning's law: $\alpha=-\dfrac{7}{6}u\beta$.
\end{itemize}
Moreover, it is natural to assume that $\beta>0$ since sediment rate increases with that of flow. 

The characteristic polynomial of $A(W)$ can be written as
\begin{equation*}
\begin{array}{ll}
p_A(\lambda)& = -\lambda[(u-\lambda)^2-gh)]+gh(\beta\lambda+\alpha)\\
& = -(\lambda -\lambda_1)(\lambda-\lambda_2)(\lambda-\lambda_3),
\end{array}
\end{equation*}
where $\lambda_1,~\lambda_2,~\lambda_3$ represent the waves speeds propagation of the current and sediment transport. The product of these
 eigenvalues, in the case $u>0$ (so $\alpha>0$), satisfies
\begin{equation*}
\lambda_1\lambda_2\lambda_3 = p_A(0)=gh\alpha <0.
\end{equation*}
This means that there exists at least a negative eigenvalue. In a sub-critical flow ($Fr<1$), one wave speed of the current is negative and
 the other one is positive. So the wave speed relating to sediment transport is positive and consequently the dune migrates downstream. Contrary,
 when the flow becomes supercritical ($Fr>1$), two waves speed of the current are positive and the wave corresponding to sediment transport
 propagates upstream so we have anti-dune.

System~\eqref{bedload system} with closed laws~(\ref{bedload Grass}-\ref{Camenen}) of $q_b$ has recently been proved to be hyperbolic within
 the range of flow data typical of practical situations~\cite{Cordier2011,Le2012}. Moreover, numerical approximations of such system are
 often based on a splitting method, solving first shallow water equations on a time step and updating afterwards the topography using the Exner
 equation. It is shown that this strategy can create spurious/unphysical oscillations. By contrast, a numerical method that solves the whole
 system at once will be called a coupled approach in the following.

This paper is organized as follows: in next section, we describe a relaxation approach of the coupled system~\eqref{bedload system} and we
 derive a numerical solver. Next, we describe the parallelization procedure we apply in order to run 2D test cases. Finally we present numerical
 test cases for moving dunes in one and two dimensions, including so-called anti-dune process.

\section{Numerical method}

In this work we want to investigate a coupled approach to numerically solve the 1D and 2D SW-Exner system. 
Such methods were already proposed in the last decade, mostly by considering extension of classical methods to non-conservative systems : Hudson
 et al. \cite{Hudson2005b} and Castro et al. \cite{castro08sediment} use a non-conservative Roe scheme based on the theory of paths, see \cite{dalmaso95definition} ; Benkhaldoun et al. \cite{Benkhaldoun09b} apply the SRNH scheme whereas Canestrelli et al. \cite{Canestrelli09} introduce the PRICE-R scheme, that is the extension to non-conservative systems of the so-called FORCE scheme. In \cite{CastroDiaz09,Benkhaldoun09,Canestrelli10}, 2D versions of the preceding works are presented. Other authors also present 2D numerical schemes based on approximate Riemann solvers as extension of HLLC scheme \cite{Murillo10,Soares10,Liang11}. In \cite{Guillard10}, the authors present an implicit procedure using the solvers presented in \cite{castro08sediment,Benkhaldoun09b}.

\medskip

Here we apply a relaxation solver to find approximate solutions of  (\ref{bedload system}). The relaxation framework for the SW-Exner model was introduced and analysed in \cite{Audusse12}. Note that a first relaxation model for SW-Exner model was introduced in \cite{Delis08} but it slightly differs from the model that we are using in this work since it is very generic and it does not use the particular form of the SW-Exner system in its definition. A relaxation solver is a particular approximate Riemann solver where the linearization is introduced in the definition of the relaxation model. By comparison to other types of approximate Riemann solvers, the relaxation framework presents some advantages since it is possible to ensure the positivity of the water height and to prove discrete energy estimates. Moreover the relaxation approach does not need a precise computation of the eigenvalues of the original system and it can be applied to conservative and non-conservative systems. The main idea if the relaxation framework is to replace the fully non-linear system (\ref{bedload system}) 
by an enlarged relaxation model that involves two types of parameters (called relaxation parameter and wave celerity parameters) 
and that satisfies the following properties
\begin{itemize}
\item The relaxation model formally tends to the original system when the relaxation parameter tends to zero.
\item The relaxation model is stable under some bounds on the wave celerity parameters.
\item The hyperbolic part of the relaxation system is linearly degenerate.
\item The related Riemann problem can be analytically solve in a (quite) easy way.
\end{itemize}
Starting from the numerical approximation at time $t^n$, the numerical procedure is then very simple : first the auxiliary quantities are computed by using the physical quantities at time $t^n$, then a homogeneous Riemann problem is solved at each interface on time step $[t^n,t^{n+1}]$ and the new physical quantities at time $t^{n+1}$ are computed by considering the mean value of the solution of these Riemann problems on each cell.\\

\medskip
The particular relaxation system that we consider in this work stands
\begin{eqnarray}
\dfrac{\partial H}{\partial t}+\dfrac{\partial Hu}{\partial x } &=& 0
\label{relaxh}\\
\dfrac{\partial Hu}{\partial t } + \dfrac{\partial}{\partial x} \left(
  Hu^2+ \Pi \right) +gH\dfrac{\partial Z}{\partial x} &=& 0
\label{relaxq}\\
\dfrac{\partial \Pi}{\partial t } + u \dfrac{\partial \Pi}{\partial x}+\dfrac{a^2}{H}\dfrac{\partial u}{\partial x} &=& \frac{1}{\lambda}(\frac{gH^2}{2}-\Pi)
\label{relaxpi}\\
\dfrac{\partial Z}{\partial t}+\dfrac{\partial \Omega}{\partial x }&=&0
\label{relaxz}\\
\dfrac{\partial \Omega}{\partial t } + \left(\dfrac{b^2}{H^2}-u^2\right)\dfrac{\partial Z}{\partial x}
+2u \dfrac{\partial \Omega}{\partial x} &=& \frac{1}{\lambda}(Q_s-\Omega)
\label{relaxqs}
\end{eqnarray}
where $\Pi$ and $\Omega$ are the auxiliary quantities (associated to fluid pressure and to sediment flux, respectively), $\lambda$ is the relaxation parameter and $a$ and $b$ are wave celerity parameters. The detailed analysis of this relaxation model and the derivation of the related relaxation scheme have been performed in \cite{Audusse12}. Let us recall some important points

\begin{itemize} 
\item The "fluid" part of the relaxation model (\ref{relaxh})-(\ref{relaxqs}), i.e. the three first equations, is nothing but the so-called Suliciu relaxation model for the classical shallow water system. The "solid" part of the relaxation model (\ref{relaxh})-(\ref{relaxqs}), i.e. the two last equations, is a classical relaxation model for a scalar conservation law and the coefficients are chossen in order to retrieve wave celerities that are centered around the water velocity $u$.
\item The stability of the relaxation model is ensured by the fact that the first order perturbation regarded parameter $\lambda$ is a diffusive perturbation. It requires some bounds on the wave celerity parameters 
\begin{equation}
a \geq H \sqrt{gH}, \qquad b \geq \sqrt{(Hu)^2+ g H^2 \partial_u Q_s}
\label{ab}
\end{equation}
\item The problem of the choice of the value of the relaxation parameter $\lambda$ is avoided by the numerical strategy that has been used for the relaxation scheme. We consider a time splitting scheme : first we only consider the right hand side of system (\ref{relaxh})-(\ref{relaxqs}) with $\lambda=0$, which only means that the auxiliary quantities instantaneously reduce to their equilibrium values ; second we compute the solution of the homogeneous (since there is no more right hand side in the equations) Riemann problem and then there is no more $\lambda$ parameter at this step.
\item The well-defined character of the relaxation scheme and the positiviy of the discrete water height may require strenghter bounds on the wave celerity parameters $a$ and $b$. 
\end{itemize}


\medskip

In \cite{Audusse12}, first numerical results show that this relaxation system is stable in some situations where the splitting approach leads to unphysical oscillations. Here we want to quantitatively analysed our results on some classical test cases but also to present the ability of the numerical model to simulate the so-called anti-dune formation. We want also to present the 2D extension of this work by using some parallelization procedure.

\section{Parallelization}

As we noted in \cite{Audusse12}, the relaxation method is very diffusive. This defect may disappear by increasing the order of the scheme. However, there is another defect thus the rise in order is not trivial. Indeed, the scheme is not well balanced. Since \cite{Bermudez94}, it is well known that the topography source term has to be treated carefully. In case of steady states at rest with varying topography, velocities might be non null. Schemes which preserve these steady states are said to be well balanced \cite{Greenberg96}. Since \cite{Bermudez94}, numerous schemes have been developed for shallow water equations \cite{Bermudez94,Bermudez98,Greenberg96,LeVeque98,Jin01a,Jin04,Jin05,Audusse04c,Audusse05,Castro07,Liang09b}... So we will have to modify our method in order to get a well-balanced scheme. But this is not trivial and unnecessary if the scheme does not catch the expected physics.

Thus the solution consists in considering finer meshes and using a parallel version of the code. Relaxation has been coded in 2D on a structured mesh. It is based on methods of line, thus parallelization is much easier than getting a high order well-balanced scheme. Parallelization has been performed thanks to domain decomposition based on the classical MPI (Message Passing Interface) library. We have adapted what is done for advection-diffusion problem in \cite{Brugeas96} to our problem.

These last decades, domain decomposition methods have increasingly attracted interest. These methods are well suited to distributed memory  parallel architectures. We decompose a domain into several sub-domains (as many sub-domains as process). In each of these sub-areas is performed calculations at the local level. The data on each side of the interfaces between sub-domains are exchanged via communication messages. The size of the domain interface is much smaller than the size of the overall problem. This domain decomposition approach leads to very natural parallelization and has the advantage to be well suited to the use of local memory.

A parallel version of the algorithm, based on the general principle of the domain decomposition adapted to structured grids, is as follows:
\begin{itemize}
 \item divide the total computational domain into as many sub-domains as process
 \item assign each sub-domain as the local domain of a process
 \item for each sub-domain determine its sub-domain neighbours
 \item iterate in time
\begin{itemize}
\item each process has to communicate with its adjacent process in order to get the flow data required for solving the equations on its local boundary cells
\item each process executes the serial code for all the computational cells lying in its local domain.
\end{itemize}
\end{itemize}

The use of MPI in case of domain decomposition is efficient. There are indeed several predefined functions that allow to create a virtual grid processes, to identify adjacent processes,... which simplify parallel programming.

\begin{figure}[ht!]
\centering
\includegraphics[angle=-0,width=0.8\linewidth]{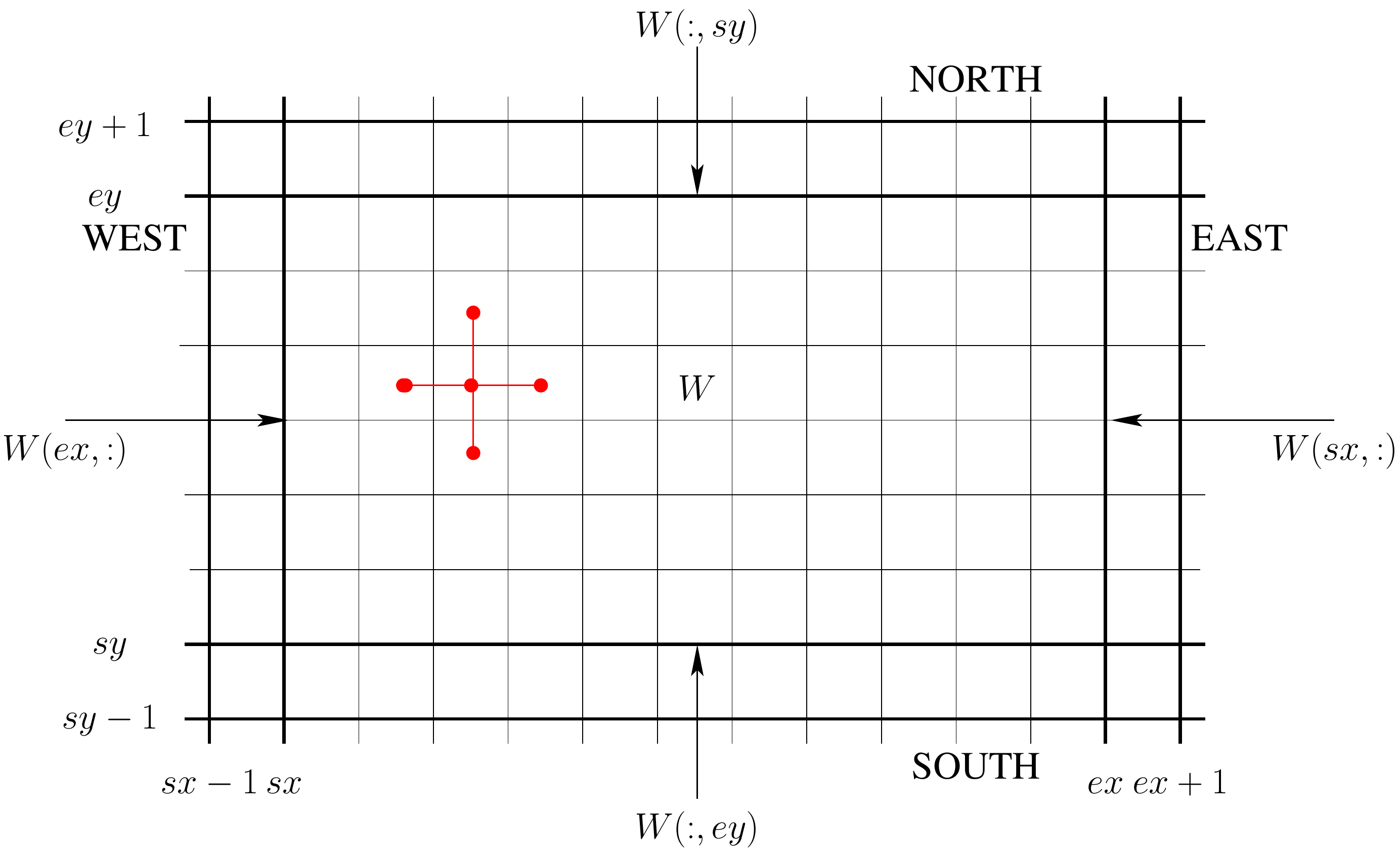}
\caption{The computational domain with its boundary conditions}
\label{fig:domain}
\end{figure}

The programming model chosen is a SPMD model (Single Program Multiple Data). One code runs on each of the sub-domains. There are as many processes as sub-domains. Each sub-domain (or process) needs to know its 8 neighbours. In the time loop, it will exchange data with 8 neighbouring interfaces and calculating within each sub-domain. We split the domain in blocks using the MPI topology feature $MPI\_CART\_CREATE$. This function returns a communicator
to which the cartesian topology information is attached. We let the function reorder the processes to match with the topology on the physical machine.
Because of the first order scheme, we only need one ghost point  to manage communications between blocks, see Figure \ref{fig:domain}. The main advantage of using the MPI function is robustness because it is well implemented in every available MPI distribution. More details concerning this domain decomposition method might be found in \cite{Brugeas96}.

In next section, the numerical method is validated on numerical tests, which is eased thanks to the code parallelization.

\section{Numerical test cases}

This numerical method has already been validated on a number of test cases \cite{Audusse12} including an analytical solution \cite{Berthon2012} integrated in the library SWASHES \cite{Delestre2013}. The idea is here to test it on more physical tests. These tests and the results will be described in what follows.

\subsection{1D dune evolution}\label{sec-1D-dune}

\begin{figure}[ht!]
\centering
\includegraphics[angle=-90,width=0.8\linewidth]{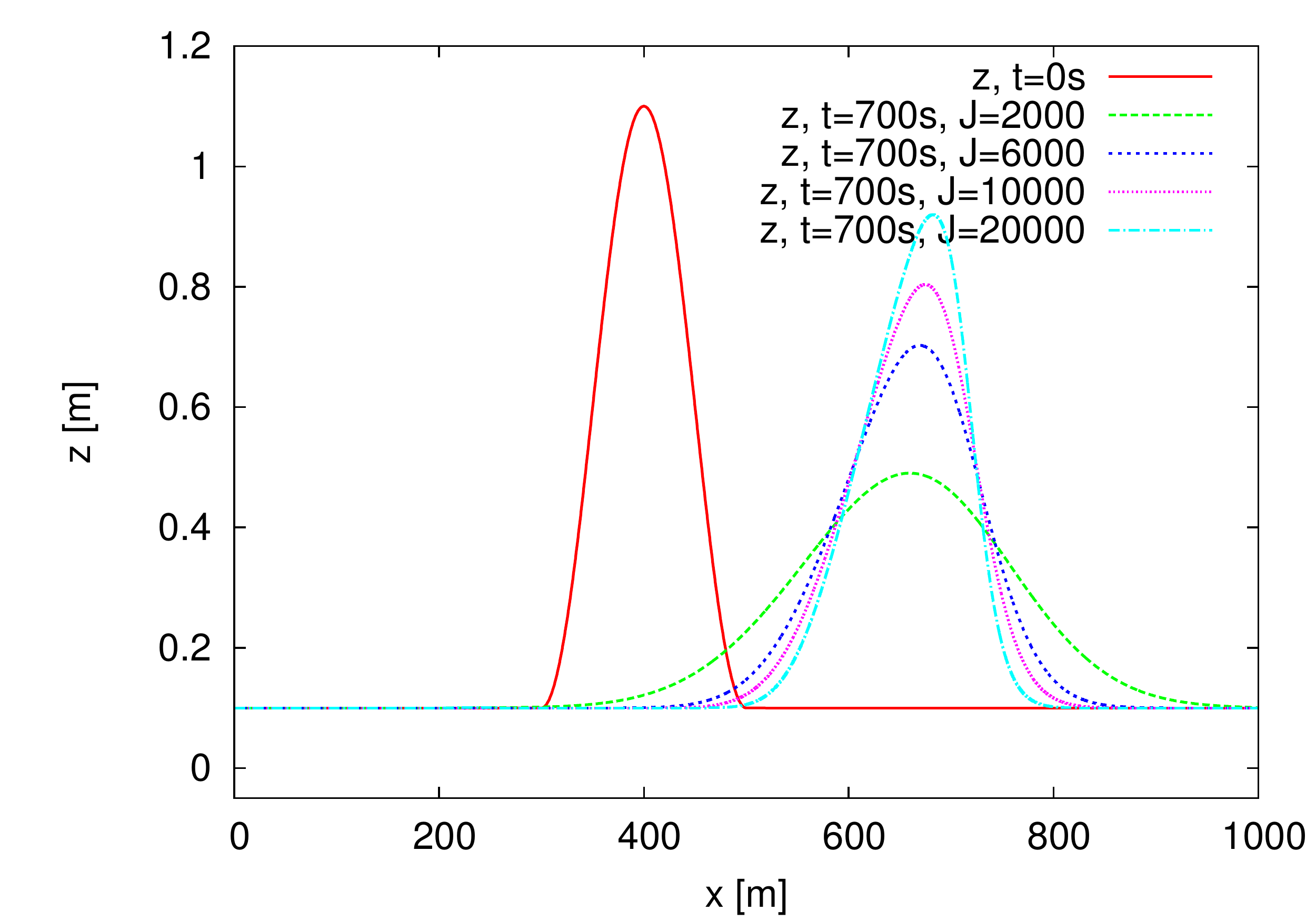}
\caption{1D dune evolution.}
\label{fig-1D-dune}
\end{figure}

We first consider a classical test case considered in several works (among others \cite{Hudson01,Hudson2003,Hudson2005b,Caleffi07,Bresch08,Delis08,castro08sediment,CastroDiaz09,Benkhaldoun09,Guillard10}), {\it i.e.} the evolution of a dune. It is a sediment transport problem in a channel of length $L=1000~\text{m}$. The initial bottom topography is a bump, we have the following initial data
\begin{equation*}
\left\{\begin{array}{l}
z_b(0,x)=\left\{\begin{array}{ll}
0.1+ \sin^2\left(\dfrac{(x-300)\pi}{200}\right) & \text{if}\;300 \leq x \leq 500~\text{m} \\
0.1 & \text{elsewhere}
\end{array}\right.\\
h(0,x)=10-z_b(0,x)\\
u(0,x)=\dfrac{q_0}{h(0,x)}
\end{array}\right.
\end{equation*}
with $q(t,0)=q_0=10~\text{m}^2/\text{s}$ the inflow discharge. Thus the flow is sub-critical. Fluid velocity increases as the bump increases  and then decreases when it decreases, so does the same $q_b$. So before the top of the bump, there is erosion ($\partial_x q_b$ is positive),  after which there is sedimentation ($\partial_x q_b$ is negative). So the bump decreases before the summit, after increases, the overall result is movement to the right. This is the evolution of the dune. For this test, we consider the Exner's law with the Grass formula  \eqref{bedload Grass}: $A_g=1$ and $m_g=3$. Thus we have a fast speed of interaction between water flow and bed-load. For the simulation  time, we have taken $T=700~\text{s}$, as in \cite{Guillard10} in order to validate our results. Four uniform grids are considered for the  discretization of the computational domain composed by $J=2000$, $6000$, $10000$ and $20000$ cells. We recover kinematics obtained in  \cite{Guillard10}. With $J=2000$ cells, we notice that the numerical method is diffusive (Fig. \ref{fig-1D-dune}). Adding cells improves the results. With these results, the 2D code parallelization is fully justified. The numerical method is suitable for sub-critical flow where  dune phenomenon occurs. With the following test, we will see that, it is suitable for supercritical flow as well.

\subsection{1D antidune evolution}\label{sec-1D-antidune}

\begin{figure}[!tb]
\centering
       \includegraphics[width=0.48\textwidth]{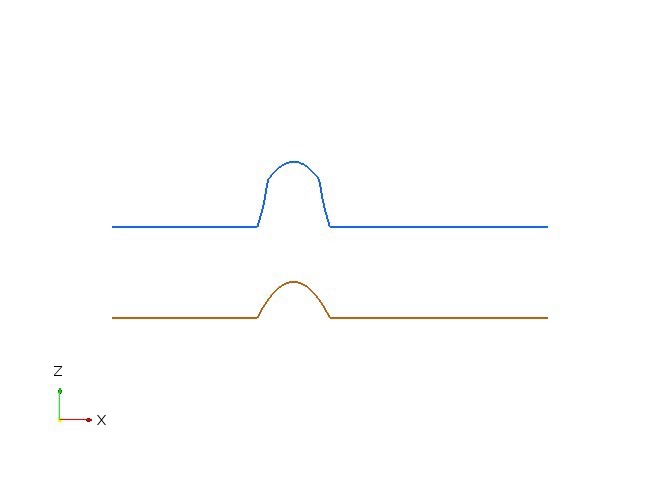}        \includegraphics[width=0.48\textwidth]{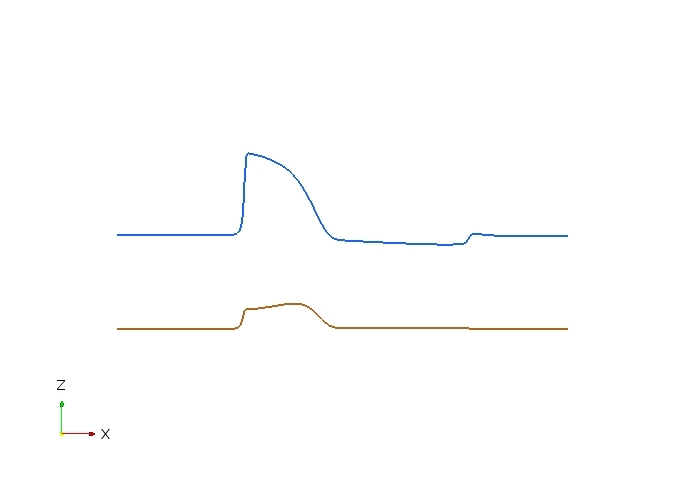}
\vspace{.25cm}
       \includegraphics[width=0.48\textwidth]{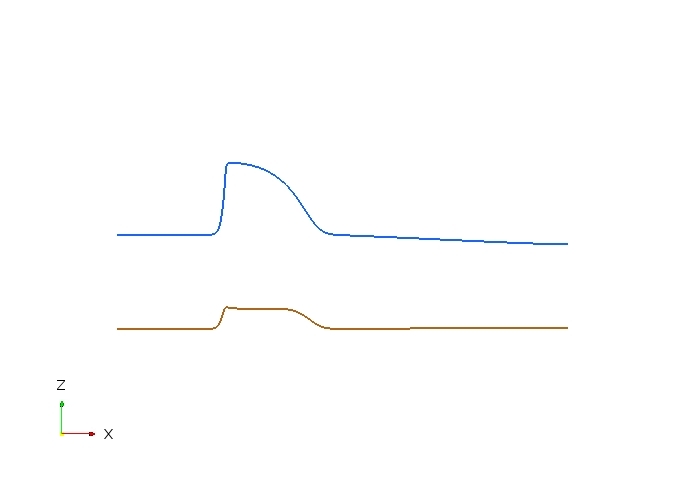}       \includegraphics[width=0.48\textwidth]{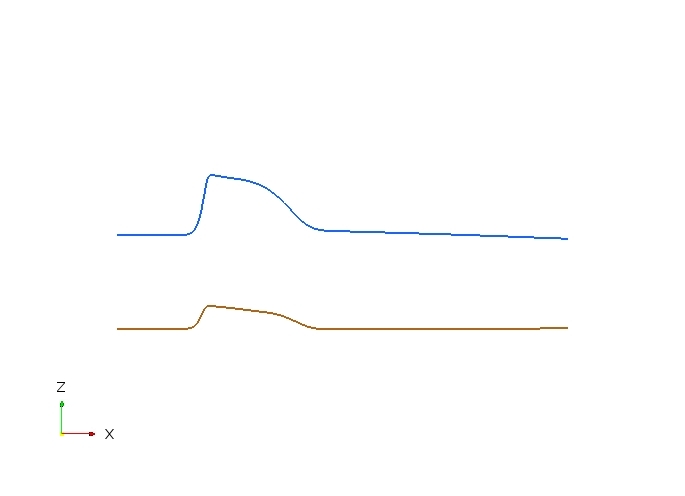}
\vspace{.25cm}
       \includegraphics[width=0.48\textwidth]{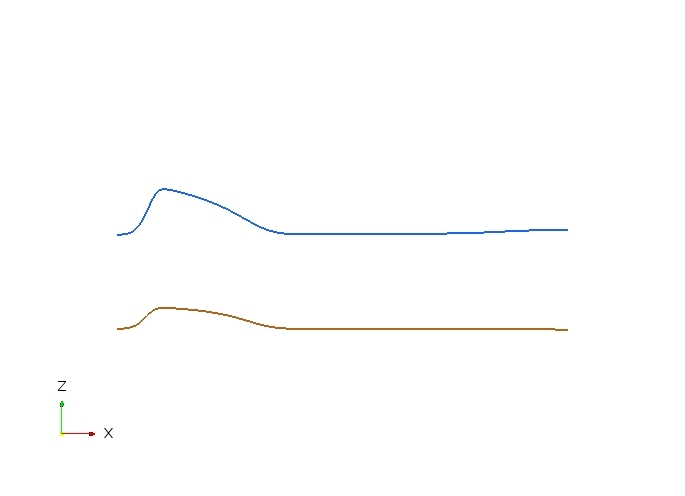}       \includegraphics[width=0.48\textwidth]{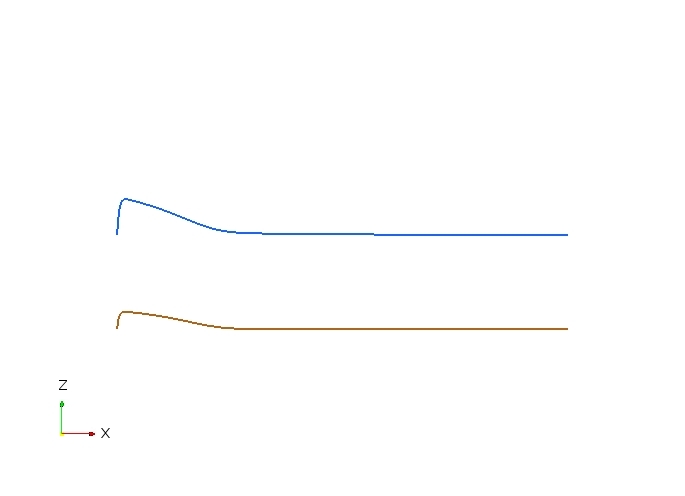}

\caption{Antidune evolution for different times 0s, 6s, 10s, 15s, 30s and 50s.}
\label{fig-antidune-1D}
\end{figure}

With this test we consider a configuration which allows the occurrence of the anti-dune phenomenon. The channel is $L=24~\text{m}$ long. The initial data are a topography with a parabolic bump
\begin{equation*}
z_b(0,x)=\left\{\begin{array}{ll}
0.2-0.05\left(x-10\right)^2  & \text{if}\;8 \leq x \leq 12~\text{m}\\
0 & \text{elsewhere}
\end{array}\right.,
\end{equation*}
a uniform discharge $q(0,x)=q_0=1.7~\text{m}^2/\text{s}$ and the water height is the stationary supercritical profile (for the Shallow Water
 equations) obtained thanks to Bernoulli's law
\begin{equation*}
\left\{\begin{array}{l}
q(t,x)=q_0\\
\dfrac{q_0^2}{2gh^2}+h+z_b=H_0
\end{array}\right.
\end{equation*}
where $H_0=\frac{q_0^2}{2gh_0^2}+h_0+z_b(0,0)$ is the total hydraulic head at inflow (Fig. \ref{fig-antidune-1D} top left). The water height at inflow is $h(t,0)=h_0=0.5~\text{m}$. This is the opposite conclusion as for the previous test: the speed decreases before the top of the bump, there is deposition, and after the summit, the speed rises again and thus causes erosion. We get anti-dune going upstream. For this test, we consider the Exner's law with the Grass formula \eqref{bedload Grass}: $A_g=0.001$ and $m_g=3$. The simulation time is $T=50~\text{s}$ and the computational domain is divided into $J=2400$ space cells. On Fig. \ref{fig-antidune-1D}, we notice the propagation of the anti-dune upstream (from right to left). This
 phenomenon has already been studied theoretically et experimently (see e.g. \cite{Kennedy1963,Parker1975,Guy1966,Huang2001,Colombini2004,Colombini2008}), but up to us, it is the first time that anti-dunes are simulated numerically.

\subsection{2D dune evolution cases}\label{sec-2D-antidune}

\begin{figure}[!htb]
\centering
       \includegraphics[width=0.5\textwidth]{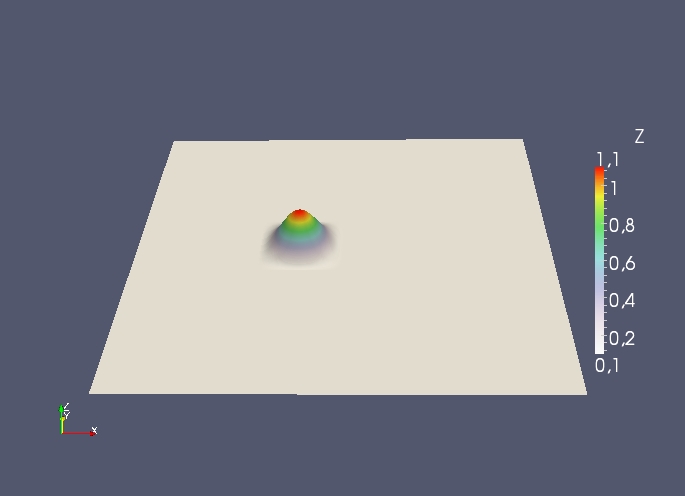}        
\caption{2D bump - initial condition.}
\label{fig.boss_2d}
\end{figure}

In these classical purely two-dimensional tests, considered in many publications (see e.g. \cite{Hudson01,Hudson05a,Delis08,Bresch08,CastroDiaz09,Benkhaldoun09,Zabsonre09,Canestrelli10,Guillard10}), we study the evolution of a bump in a channel (Fig. \ref{fig.boss_2d}). The dimensions of the domain are $(L_x\times L_y)=(1000\text{m}\times 1000\text{m})$. The initial data common to each test are
\begin{equation*}
 \left\{\begin{array}{l}
z_b(0,x,y)=\left\{\begin{array}{ll}
                           0.1+\sin \left(\pi \dfrac{(x-300)}{200}\right)^2 \sin\left(\pi\dfrac{(y-400)}{200}\right)^2 & \text{if}\,
 (x,y)\in[300,500\text{m}]\times [400,600\text{m}]\\
        	  0.1 & \text{elsewhere}
                          \end{array}\right.\\
h(0,x,y)=10-z_b(0,x,y)\\
u(0,x,y)=\dfrac{q_0}{h(0,x,y)}\\
v(0,x,y)=0         
 \end{array} \right.
\end{equation*}
with $u$ (resp. $v$) the velocity in $x$ (resp. $y$) direction and $q_0=10\text{m}^2/\text{s}$ the upstream constant discharge in $x$ direction. In all other boundaries, we assume free flow conditions. For these tests, we consider the Exner's law with the 2D Grass formula with $m_g=3$ and two different values: $A_g=1$ and then $A_g=0.1$. Thus we get a fast and an intermediate ground/flow interaction.

\subsubsection{Fast interaction -- $A_g=1$}

\begin{figure}[!htb]
\begin{center}
       \includegraphics[width=0.48\textwidth]{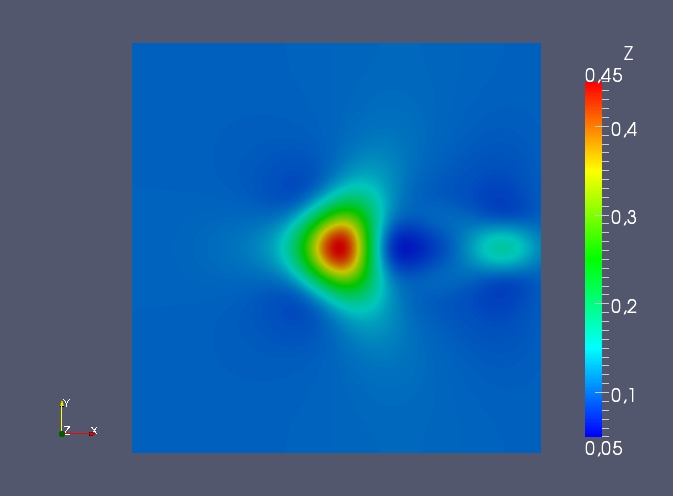}        \includegraphics[width=0.48\textwidth]{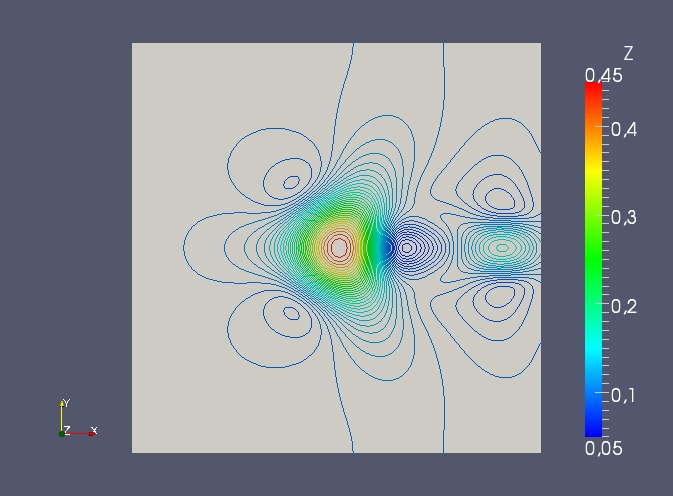}

\caption{2D bump at $T=500~\text{s}$ for $A_g=1$.}
\label{fig-bump-2D-As1}
\end{center}
\end{figure}

For this test, intial data are those described previously with $A_g=1$. The simulation time is $T=500~\text{s}$ and the computational domain is divided into $(J\times K)=(4000\times 4000)$ cells. At the end of the simulation (Fig. \ref{fig-bump-2D-As1}), we recover the shape obtained in \cite{CastroDiaz09} and \cite{Guillard10} for this value of $A_g$.

\subsubsection{Intermediate interaction -- $A_g=0.1$}

\begin{figure}[!htb]
\begin{center}
       \includegraphics[width=0.48\textwidth]{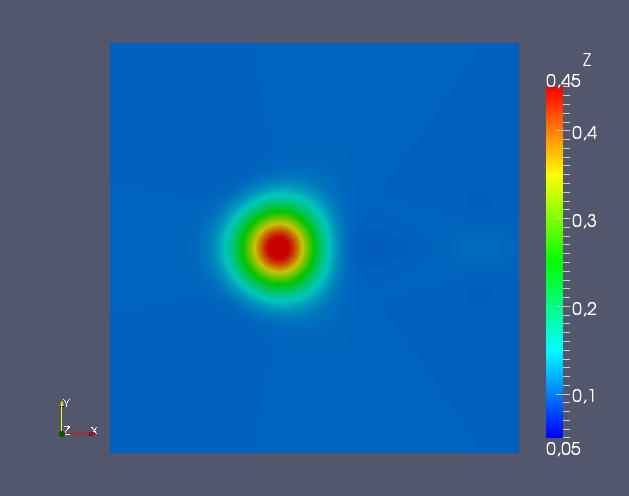}        \includegraphics[width=0.48\textwidth]{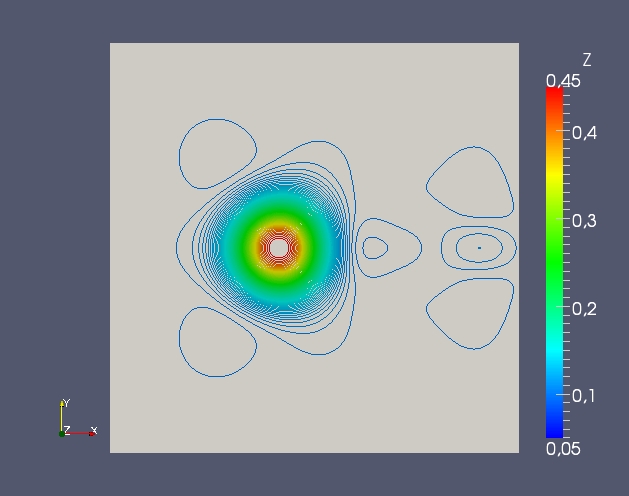}
\vspace{.25cm}
       \includegraphics[width=0.48\textwidth]{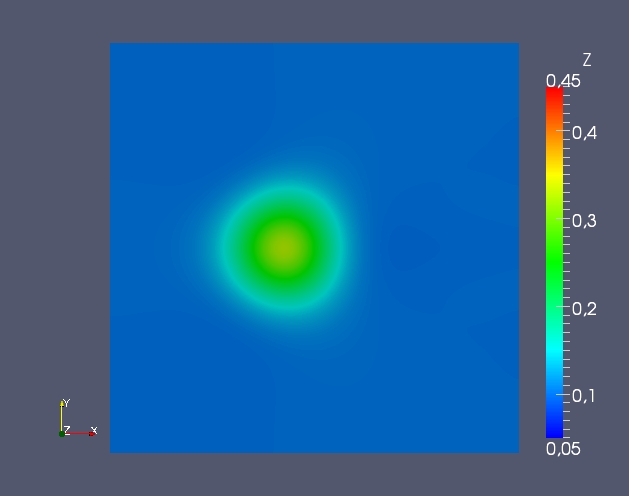}       \includegraphics[width=0.48\textwidth]{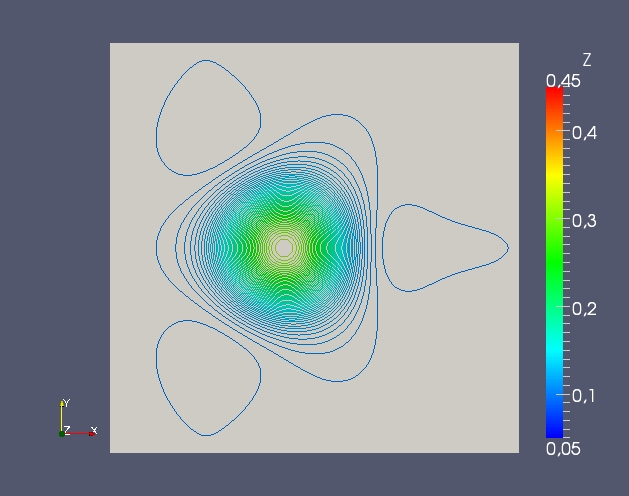}

\caption{2D bump evolution at $T=500~\text{s}$ (up) and $T=1000~\text{s}$ (down) for $A_g=0.1$.}
\label{fig-bump-2D-As01}
\end{center}
\end{figure}

For this test, intial data are those described previously with $A_g=0.1$. The simulation time is $T=1000~\text{s}$ and the computational domain is divided into $(J\times K)=(4000\times 4000)$ cells. At the location of the of the bump, we notice that the method is diffusive (Fig. \ref{fig-bump-2D-As01}): indeed, we have circles on the pictures representing contours. However, as for the previous test, the overall shape is correct.

\begin{figure}[!htb]
\begin{center}
       \includegraphics[width=0.48\textwidth]{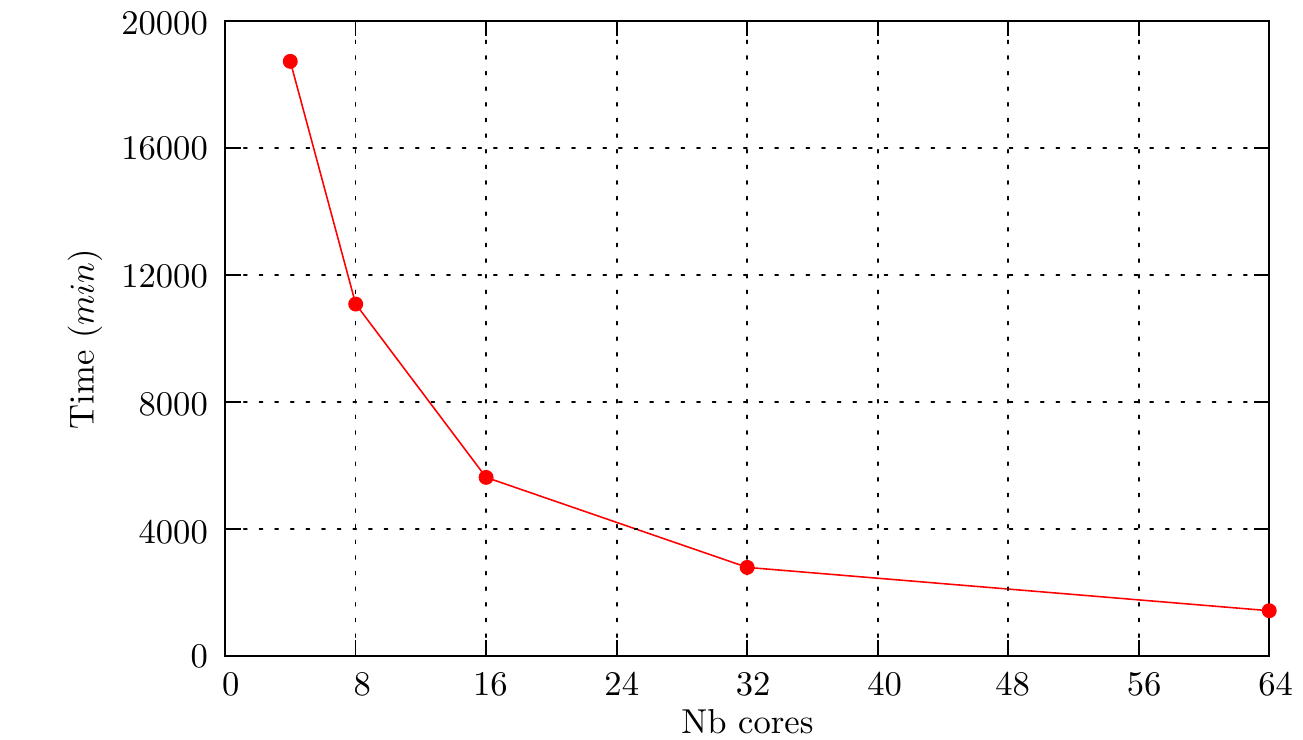}        \includegraphics[width=0.48\textwidth]{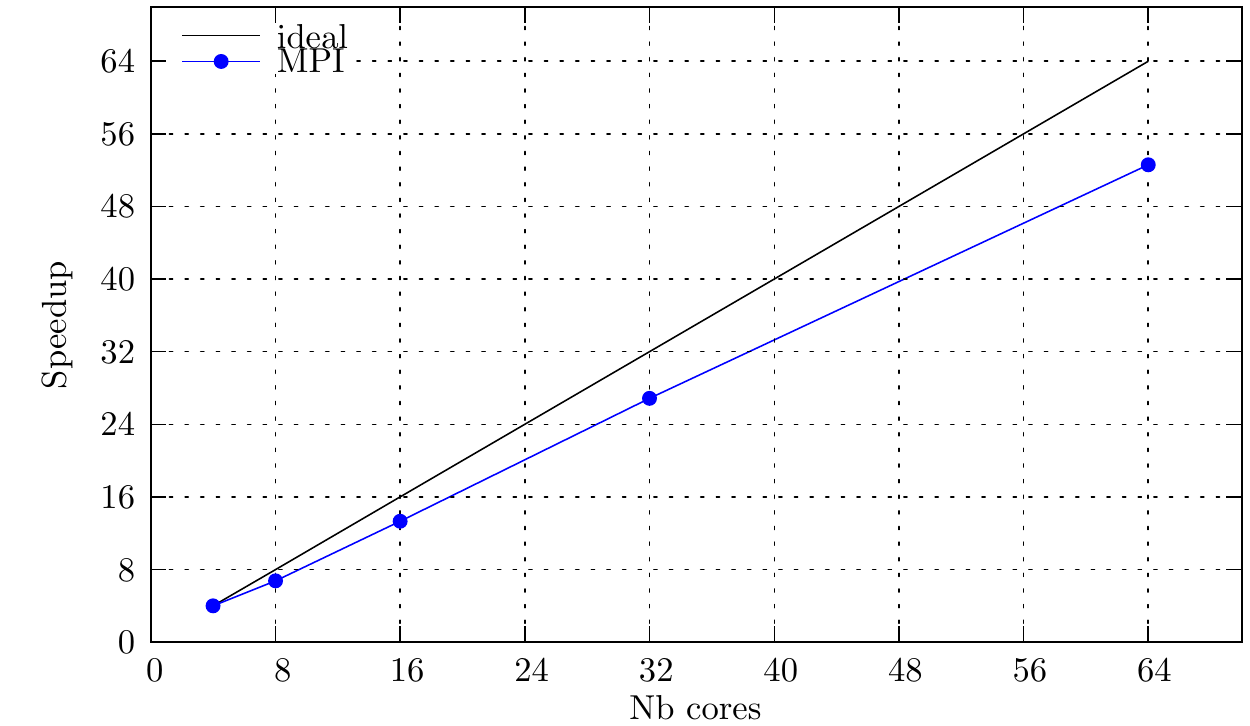}

\caption{CPU time and speedup for $A_g=0.1$.}
\label{fig-performance}
\end{center}
\end{figure}

In order to evaluate the speed-up, we ran the simulation on $4$, $8$, $16$, $32$ and $64$ cores. The computation time and acceleration are represented on Fig. \ref{fig-performance}. Usually, speed-up is defined by $S=T_1/T_p$, with $p$ the number of processors, $T_1$ the execution time of the sequential algorithm and $T_p$ the execution time of the parallel algorithm with $p$ processors. Thus speed-up is a very classical {\it criteria} to know how much a parallel algorithm is faster than a corresponding sequential algorithm (see among others \cite{Karniadakis03,ElRewini05} and \cite{Pacheco11}). A linear speed-up or ideal speed-up is obtained when $S_p=p$. In our case, performing the sequential code was too consuming, thus we have used $T_4$ as the reference. The "speed-up" is calculated as what follows: $S=4T_4/T_p$. We get a speed-up not far from the ideal one. For future improvements and to take advantage of the power of massively parallel machines of the last generation, hybrid parallelization should be considered. Using OpenMP with MPI could increase benefits like: memory saving, better load balancing and better adequacy to the hardware specificities.

\section{Conclusions and perspectives}

The objective of this work was to prove that the relaxation solver that we introduced in \cite{Audusse12} is adapted to reproduce 
\begin{itemize}
\item well known physical phenomena as the anti-dune process that are not so easy to handle at the numerical level, in particular when soft coupling ({\it{i.e.}} independent solvers or fluid and solid parts) are considered as it is mostly the case for industrial applications,
\item two dimensional situations when coupled with a parallelization procedure.
\end{itemize}

For future works, some improvements have to be performed at the numerical level : high order, well balancing... 
But the most promising way seems to propose some improvements for the SWE model \eqref{bedload system} itself. We see at least two directions that can be investigated. The first one is to introduce more physics in the sediment layer by considering a more complete system for which energy equation can be exhibited. This can be done by considering formal reduction procedure starting from three dimensional model (in the spirit of \cite{Gerbeau01} for classical shallow water system) or by directly introducing adapted multilayer models, see for example \cite{spinewine08}.
The other one is to improve the fluid model.
Indeed, performing linear stability analysis can show that the Saint-Venant equations, which is a averaged models, do not admit instability rising on the bed from a suitable initial perturbation contrary to those observed in real life. This relative failure of the models led to the consideration of a full fluid flow model, in which, rather than supposing that the flow is shear free and that viscous effects were confined to a turbulent boundary layer, rotational effects were considered, and a model of turbulent shear flow incorporating an eddy viscosity, together with the Exner equation for bedload transport, was adopted (see e.g. {\cite{Colombini2004,Lagree2003,Devauchelle2007,Devauchelle2010,Charru2002,Kouakou2005,Camporeale2009}). 

\section*{Acknowledgments}
O. Delestre, R. Serra, H. Le Minh thank CaSciModOT federation (Calcul Scientifique et Mod\'elisation Orl\'eans Tours) and AMIES agency (Agence
 pour les Math\'ematiques en Interaction avec l'Entreprise et la Soci\'et\'e) for financial support. E. Audusse and M. Masson-Fauchier thank
 PEPS project (Projets Exploratoires Premier Soutien pour les Interactions Math\'ematiques-Informatique-Ing\'enierie) initiated by INSMI for
 financial support.


\bibliographystyle{plain}

\end{document}